\documentclass[11pt, a4paper]{amsart}
\usepackage{amsmath}
\usepackage{amssymb}
\usepackage{color}

\newtheorem{teo}{Theorem}[section]

\theoremstyle{definition}

\def\a{\alpha}

\def\R{\mathbb R}

\begin{document}
\title[Asymptotics for inhomogenous equations with memory. Small dimensions]{Decay/growth rates for inhomogeneous heat equations with memory. The case of small dimensions}

\author[C. Cort\'{a}zar,  F. Quir\'{o}s \and N. Wolanski]{Carmen Cort\'{a}zar,  Fernando  Quir\'{o}s, \and Noem\'{\i} Wolanski}

\address{Carmen Cort\'{a}zar\hfill\break\indent
Departamento  de Matem\'{a}tica, Pontificia Universidad Cat\'{o}lica
de Chile \hfill\break\indent Santiago, Chile.} \email{{\tt
ccortaza@mat.puc.cl} }

\address{Fernando Quir\'{o}s\hfill\break\indent
Departamento  de Matem\'{a}ticas, Universidad Aut\'{o}noma de Madrid,
\hfill\break\indent 28049-Madrid, Spain,
\hfill\break\indent and Instituto de Ciencias Matem\'aticas ICMAT (CSIC-UAM-UCM-UC3M),
\hfill\break\indent 28049-Madrid, Spain.} \email{{\tt
fernando.quiros@uam.es} }

\address{Noem\'{\i} Wolanski \hfill\break\indent
IMAS-UBA-CONICET, \hfill\break\indent Ciudad Universitaria, Pab. I,\hfill\break\indent
(1428) Buenos Aires, Argentina.} \email{{\tt wolanski@dm.uba.ar} }



\date{}

\begin{abstract}
We study the decay/growth rates in all  $L^p$ norms of solutions to  an inhomogeneous nonlocal heat equation in $\mathbb{R}^N$ involving a Caputo  $\alpha$-time derivative and a power $\beta$ of the Laplacian when the spatial dimension is small, $1\le N\le 4\beta$, thus completing the already available results for large spatial dimensions. Rates depend not only on $p$, but also on the space-time scale and on the time behavior of the spatial $L^1$ norm  of the forcing term.
\end{abstract}

\maketitle

\section{Introduction}
 \setcounter{equation}{0}

\subsection{Goal}  This paper is devoted to study the large-time behavior of the $L^p$-norms, $p\in[1,\infty]$,  of the  solution to the \emph{fully nonlocal} problem
 \begin{equation}\label{eq-f}
 	\partial_t^\a u+(-\Delta)^\beta u=f\quad\mbox{in }Q:=\R^N\times(0,\infty),\qquad
 	u(\cdot,0)=0\quad\mbox{in }\R^N,
\end{equation}
in the case of small dimensions, $1\le N\le 4\beta$. It may be regarded as a companion to~\cite{Cortazar-Quiros-Wolanski-2022}, where the case  of large dimensions, $N> 4\beta$, was treated.  Here, $\partial_t^\alpha$, $\alpha\in(0,1)$, denotes the so-called Caputo $\alpha$-derivative, introduced independently by many authors using different points of view, see for instance~\cite{Caputo-1967,Dzherbashyan-Nersesian-1968,Gerasimov-1948,Gross-1947,Liouville-1832,Rabotnov-1966}, which is defined for smooth functions by
$$
\displaystyle\partial_t^\alpha u(x,t)=\frac1{\Gamma(1-\alpha)}\,\partial_t\int_0^t\frac{u(x,\tau)- u(x,0)}{(t-\tau)^{\alpha}}\, {\rm d}\tau,
$$
and $(-\Delta)^\beta$, $\beta\in(0,1]$,  is the usual $\beta$ power of the Laplacian, defined for smooth functions by
$(-\Delta)^{s}=\mathcal{F}^{-1}(|\cdot|^{2s}\mathcal{F})$, where $\mathcal{F}$ denotes Fourier transform; see for instance~\cite{Stein-book-1970}.
Equations of this kind,  nonlocal both in space and time, have been proposed recently to model situations with long-range interactions and memory effects; see for example~\cite{Cartea-delCastilloNegrete-2007,Compte-Caceres-1998,delCastilloNegrete-Carreras-Lynch-2004,delCastilloNegrete-Carreras-Lynch-2005,Metzler-Klafter-2000,Zaslavsky-2002}.

Problem~\eqref{eq-f} does not have in general a classical solution, unless the forcing term $f$ is required to satisfy certain smoothness assumptions. However, if $f\in L^\infty_{\rm loc}([0,\infty):L^1(\mathbb{R}^N))$,  it has a solution in a generalized sense,
defined by Duhamel's type formula
\begin{equation}\label{eq-formula}
 		u(x,t)=\int_0^t\int_{\R^N}Y(x-y,t-s)f(y,s)\,{\rm d}y{\rm d}s,
\end{equation}
with $Y=\partial_t^{1-\a} Z$, where $Z$ is the fundamental solution of the Cauchy problem; see \cite{Eidelman-Kochubei-2004,Kemppainen-Siljander-Zacher-2017}.
Throughout the paper we assume the integral size condition
 \begin{equation}\label{eq-hypothesis f}
 \|f(\cdot,t)\|_{L^1(\mathbb{R}^N)}\le \frac C{(1+t)^\gamma}\quad\text{for some }\gamma\in\mathbb{R},
 \end{equation}
 and  deal with solutions of this kind, denoted in the literature as \emph{mild} solutions~\cite{Kemppainen-Siljander-Zacher-2017,Pruss-book}.

\medskip

\noindent\emph{Notation. } As is common in asymptotic analysis, $g\asymp h$ will mean that there are constants $\nu,\mu>0$ such that
$\nu h\le g\le \mu h$.

 \subsection{The kernel $Y$. Critical exponent}\label{subsect-estimates W}

Our study is based in a careful analysis of the singular integral~\eqref{eq-formula} defining the mild solution. Hence, having good estimates for the kernel $Y$ is essential. Such estimates are fortunately available in~\cite{Kemppainen-Siljander-Zacher-2017}, and are recalled next.

The kernel $Y$ has a self-similar form,
\begin{equation}
 \label{eq:Y.selfsimilar}
 Y(x,t)=t^{-\sigma_*}G(\xi),\quad \xi=x t^{-\theta},\qquad \sigma_*:=1-\alpha+N\theta,\quad \theta:=\frac{\alpha}{2\beta }.
\end{equation}
Its profile $G$ is positive, radially symmetric, and smooth outside the origin, and if $N\le 4\beta$ satisfies close to the origin the sharp estimates
  \begin{equation}\label{eq:estimates.G.origin}
 	G(\xi)\asymp \begin{cases}1,&N<4\beta,\\[-0.3cm]
&\hskip2.5cm \beta\in (0,1],\quad |\xi|\le1,\\[-0.3cm]
 		\big(1+\big|\log|\xi|\big|\big),&N=4\beta,
 \end{cases}
 \end{equation}
It is here that we find the main difference with respect to the case of high dimensions, $N>4\beta$, for which
$$
G(\xi)\asymp C |\xi|^{4\beta-N},\qquad \beta\in(0,1],\quad |\xi|\le 1.
$$
Away from the origin we have the same (sharp) behavior as in high dimensions,
 \begin{equation}\label{eq:exterior.estimate.G}
 	G(\xi)\asymp \begin{cases} |\xi|^{-(N+2\beta)},&\beta\in(0,1) \\[-0.3cm]
&\hskip2.5cm   |\xi|\ge1.\\[-0.3cm]
 |\xi|^{\frac{(N-2)(\alpha-1)}{(2-\alpha)}}{\rm e}^{-\sigma|\xi|^{\frac2{2-\a}}},&\beta=1,
  	\end{cases}
\end{equation}
In particular,   since $|\xi|^{\frac{(N-2)(\alpha-1)}{(2-\alpha)}}\exp(-\sigma|\xi|^{\frac2{2-\a}})\le C_\nu|\xi|^{-(N+2\beta)}$ if $|\xi|\ge\nu$,  we have the exterior bound
\begin{equation}
  \label{eq:exterior.estimate.Y}
0\le Y(x,t)\le C_\nu t^{2\alpha-1}|x|^{-(N+2\beta)}\quad\text{if } |x|\ge \nu t^\theta,\ t>0,\quad \beta\in(0,1].
\end{equation}

Notice that $Y(\cdot,t)\in L^p(\mathbb{R}^N)$ for all $p\in [1,\infty]$ if $N<4\beta$ and for all $p\in[1,\infty)$ if $N=4\beta$. Moreover, in all these cases we have
\begin{equation}\label{eq:p.norm.Y}
\|Y(\cdot,t)\|_{L^p(\mathbb{R}^N)}=C t^{-\sigma(p)}\quad\text{for all }t>0, \qquad\text{where } \sigma(p):=\sigma_*-\frac{N\theta}p,\quad p\in[1,\infty].
\end{equation}
We remark that $\sigma(p)<1$, and hence $Y\in L_{\rm loc}^1([0,\infty):L^p(\R^N))$, if and only if  $p$ is \emph{subcritical},
\begin{equation}
\label{eq:subcritical.range}
\tag{S}
p\in [1,\infty]\text{ if }N<2\beta,\quad p\in[1,p_c)\text{ if }N\ge 2\beta, \quad \text{where }
p_c=\begin{cases}
N/(N-2\beta)&\text{if }N>2\beta,\\
\infty&\text{if }N=2\beta.
\end{cases}
\end{equation}
If $p$ is not subcritical we need some extra assumption on  the forcing term to guarantee that $u(\cdot,t)\in L^p(\mathbb{R}^N)$. In the present paper we will use two different such extra hypotheses,  the pointwise condition
\begin{equation}
\label{eq:decay.condition}
|f(x,t)|\le C|x|^{-N}(1+t)^{-\gamma}\quad \text{for }|x|\text{ large},
\end{equation}
and the integral condition
\begin{equation}
\label{eq:q.integrability.f}
\|f(\cdot,t)\|_{L^q(\R^N)}\le C(1+t)^{-\gamma}\quad\text{for some }q>q_{\rm c}(p):=\frac{N}{2\beta +\frac Np}.
\end{equation}
Note that $q_{\rm c}(p_{\rm c}) =1$.
We do not claim that these conditions are optimal; but they are not too restrictive, and are easy enough to keep the proofs simple.

 \subsection{Precedents and statement of results}\label{subsect-previous results}
As a first precedent we have~\cite{Kemppainen-Siljander-Zacher-2017}, where  the authors study the problem in the \emph{integrable in time} case $\gamma>1$ and  prove, when $p$ is subcritical,   that
\begin{equation*}
\label{eq:result.KSZ.subcritical}
\lim_{t\to\infty}t^{\sigma(p)}\|u(\cdot,t)-M_\infty Y(\cdot,t)\|_{L^p(\mathbb{R}^N)}=0,\quad
\text{where }M_\infty:=\int_0^\infty\int_{\mathbb{R}^N} f(x,t)\,{\rm d}x{\rm d}t<\infty,
\end{equation*}
which, as long as $M_\infty\neq0$, implies in particular the sharp decay rate
\begin{equation}
\label{eq:subcritical.global.rate.precedents} 	
\|u(\cdot,t)\|_{L^p(\R^N)}\asymp t^{-\sigma(p)}.
\end{equation}

A second and very recent precedent is~\cite{Cortazar-Quiros-Wolanski-2022}, where we study the case of large dimensions, $N>4\beta$. In sharp contrast with the case in which the time derivative is local, and due to the effect of memory, we proved there that the decay/growth rates are not the same in different space-time scales. We already found this phenomenon for the homogeneous Cauchy problem (with nontrivial initial datum) in~\cite{Cortazar-Quiros-Wolanski-2021a,Cortazar-Quiros-Wolanski-2021b}.
In the present paper we have to face the same difficulty, which forces us to study separately the rates in exterior regions, $|x|\ge\nu t^\theta$ with $\nu>0$, compact sets or intermediate  regions $ |x|\asymp g(t)$ with $g(t)\to\infty$ and $g(t)=o(t^\theta)$. The main difference in comparison with the case of large dimensions is that now we have to deal with additional critical behaviors associated with the critical dimensions $N=2\beta$ and $N=4\beta$.

Our results are summarized as follows:
\begin{teo}[\sc Exterior regions]
\label{teo-exterior}
Let $1\le N\le 4\beta$.
Let $f$ satisfy \eqref{eq-hypothesis f} and also \eqref{eq:decay.condition} if $p$ is not subcritical.  Let $u$ be the mild solution to \eqref{eq-f}. For all $\nu>0$ there is a constant $C$ such that
	\begin{equation*}\label{eq-bound exterior}
	\|u(\cdot,t)\|_{L^p(\{|x|\ge \nu t^{\theta}\})}\le C \begin{cases}
	t^{-\sigma(p)+1-\gamma},&\gamma<1,\\
t^{-\sigma(p)}\log t,&\gamma=1,\\
t^{-\sigma(p)},&\gamma>1.
\end{cases}
\end{equation*}
These estimates are sharp.
\end{teo}
When we say that the estimates are sharp we mean that there are  forcing functions $f$ such that the estimates from below hold.

\noindent\emph{Remark. } For $p\in[1,p_{\rm c})$ and $\gamma>1$ the result follows from~\eqref{eq:subcritical.global.rate.precedents}.

It is in compact regions where we find critical behaviors associated with the dimension.

\begin{teo}[\sc Compact sets]\label{teo-compactos}
Let  $f$ satisfy \eqref{eq-hypothesis f}. If $p$ is not subcritical, assume  also  \eqref{eq:q.integrability.f}  with $\gamma$ as in \eqref{eq-hypothesis f}. Let $u$ be the mild solution to \eqref{eq-f}. For every compact set $K$ there exists a constant $C$ such that
$$
\|u(\cdot,t)\|_{L^p(K)}\le C \begin{cases}
t^{-\sigma_*+1-\gamma},&\gamma< 1,\\
	t^{-\sigma_*}\log t,&\gamma= 1,\hskip2.4cm N<2\beta,\\
	t^{-\sigma_*},&\gamma> 1,\\[0.3cm]
t^{-\gamma}\log t,&\gamma\le\sigma_*= 1,\\[-0.3cm]
&\hskip3.5cm  N=2\beta,\\[-0.3cm]
	t^{-\sigma_*},&\gamma>\sigma_*= 1,\\[0.3cm]
t^{-\gamma},&\gamma< \sigma_*,\\[-0.3cm]
&\hskip3.5cm  2\beta<N<4\beta,\\[-0.3cm]
	t^{-\sigma_*},&\gamma\ge \sigma_*,\\[0.3cm]
t^{-\gamma},&\gamma< \sigma_*=1+\alpha,\\[-0.3cm]
&\hskip3.5cm  N=4\beta.\\[-0.3cm]
t^{-\sigma_*}\log t,&\gamma\ge \sigma_*=1+\alpha,
\end{cases}
$$
These estimates are sharp.
\end{teo}

\noindent\emph{Remark. } Though the decay rates in compact sets are the same in all $L^p$ norms, the value of $p$ still plays a role in our proofs, since we will need an extra decay assumption when $p$ is not subcritical.

As expected, the rates in intermediate regions, between compact sets and exterior regions, are intermediate between the ones in such scales, and hence will also exhibit critical behaviors associated with the dimension.

\begin{teo}[\sc  Intermediate regions]
\label{teo-intermediate}
Let $f$ satisfy \eqref{eq-hypothesis f} and also \eqref{eq:decay.condition} if $p$ is not subcritical. Let $g(t)\to\infty$ be such that $g(t)=o(t^{\theta})$. Let $u$ be the mild solution to \eqref{eq-f}. For all $0<\nu<\mu<\infty$ there exists a constant $C$ such that
$$
\begin{aligned}
\|u&(\cdot,t)\|_{L^p(\{\nu <|x|/g(t)<\mu\})}\\[0.3cm]
&\le C g(t)^{N/p}\begin{cases}
	t^{-\sigma_*+1-\gamma},&\gamma<1,\\
	t^{-\sigma_*}\log t,&\gamma=1, \hskip.3cm  N<2\beta,\\
	t^{-\sigma_*},&\gamma>1,\\[0.3cm]
	\log\big(t^\theta/g(t)\big)t^{-\gamma},&\gamma<1,\\
	t^{-\sigma_*}\log t,&\gamma=1, \hskip.3cm N=2\beta,\\
	t^{-\sigma_*},&\gamma>1,\\[0.3cm]
	g(t)^{(1-\sigma_*)/\theta}t^{-\gamma},&\gamma<1,\\
	\max\big\{g(t)^{(1-\sigma_*)/\theta}t^{-1},t^{-\sigma_*}\log t\big\},&\gamma=1,  \hskip.3cm 2\beta<N<4\beta,\\
	\max\big\{g(t)^{(1-\sigma_*)/\theta}t^{-\gamma},t^{-\sigma_*}\big\},&\gamma>1,\\[0.3cm]
	\max\big\{g(t)^{(1-\sigma_*)/\theta}t^{-\gamma},
	\log\big(t^\theta/g(t)\big)t^{-\sigma_*+1-\gamma}\big\},
	&\gamma<1,\\
	\max\big\{g(t)^{(1-\sigma_*)/\theta}t^{-1},\log\big(t^\theta/g(t)\big)
	t^{-\sigma_*}\log t\big\},&\gamma=1,\hskip.3cm N=4\beta.\\
	\max\big\{g(t)^{(1-\sigma_*)/\theta}t^{-\gamma},
	\log\big(t^\theta/g(t)\big)t^{-\sigma_*}\big\},&\gamma>1, 	
	\end{cases}
\end{aligned}
$$
These estimates are sharp.
\end{teo}

We also obtain results that connect the behavior in compact sets and exterior
regions thus getting the decay rate in $L^p(\R^N)$.
\begin{teo}[\sc Global results]
\label{teo-conexion}
Assume \eqref{eq-hypothesis f},  and also \eqref{eq:decay.condition}  and
\eqref{eq:q.integrability.f}
with $\gamma$ as in \eqref{eq-hypothesis f} if $p$ is not subcritical. Let $u$ be the mild solution to \eqref{eq-f}.

\noindent{\rm (i)} If $p$ is subcritical, there is a constant $C$ such that
$$
\|u(\cdot,t)\|_{L^p(\R^N)}\le C \begin{cases}
t^{-\sigma(p)+1-\gamma},&\gamma<1,\\
t^{-\sigma(p)}\log t,&\gamma=1,\\
t^{-\sigma(p)},&\gamma>1.
\end{cases}
$$

\noindent{\rm (ii)} If $p$ is critical, $N\ge 2\beta$, $p=p_{\rm c}$, there is a constant $C$ such that
$$
\|u(\cdot,t)\|_{L^p(\R^N)}
\le C \begin{cases}
t^{-\gamma}\log t,&\gamma\le1, \\
t^{-1},&\gamma>1.
\end{cases}
$$
\noindent{\rm (iii)} If $p$ is supercritical, $N>2\beta$, $p>p_{\rm c}$, there is a constant $C$ such that
$$
\|u(\cdot,t)\|_{L^p(\R^N)}
\le C \begin{cases}
t^{-\gamma},&\gamma<\sigma(p),\\	  	      	
t^{-\sigma(p)},&\gamma\ge\sigma(p),\hskip1.48cm p<\infty \text{ or } N<4\beta,\\
t^{-\sigma_*}\log t,&\gamma\ge \sigma_*=1+\alpha,\quad  p=\infty,\ N=4\beta.
\end{cases}
$$
These estimates are sharp.	
\end{teo}

Thus, if $p$ is subcritical or if $p$ is supercritical and $\gamma$ is large, the global decay rate is determined by the behavior in exterior domains, while if $p$ is supercritical and $\gamma$ is small it is dictated by the behavior in compact sets. In the critical case with $\gamma$ large, compact sets rule the game. The same is true for $\gamma$ small in the critical dimension $N=2\beta$. However, if $\gamma$ is small and $2\beta<N\le 4\beta$, then the whole intermediate region plays a role in the determination of the global rate.

Once sharp decay/growth rates are available, the next natural step is to obtain the asymptotic profile, in each space-time scale, after correcting the solution with the rate. Such goal is the content of the forthcoming paper~\cite{Cortazar-Quiros-Wolanski-2021-Preprint}.

\section{Exterior region}\label{sect-exterior}
 \setcounter{equation}{0}

This section is devoted to the proof of Theorem~\ref{teo-exterior}, which gives the decay/growth rates in all $L^p$ norms of the mild solution $u$ to \eqref{eq-f} in exterior regions, $\{(x,t)\in Q:|x|\ge \nu t^\theta\}$, $\nu>0$.

\begin{proof}[Proof of Theorem~\ref{teo-exterior}.] The starting point will be always Duhamel's type formula~\eqref{eq-formula}.

If $p$ is subcritical, that is, if \eqref{eq:subcritical.range} holds, using the size condition  \eqref{eq-hypothesis f} and~\eqref{eq:p.norm.Y} we get the global estimate
\begin{equation}
\label{eq:global.subcritical}
\begin{aligned}
	\|u(\cdot,t)\|_{L^p(\R^N)}&\le \int_0^t\|Y(\cdot,t-s)\|_{L^p(\R^N)}\|f(\cdot,s)\|_{L^1(\R^N)}\,{\rm d}s\le C\int_0^t(t-s)^{-\sigma(p)}(1+s)^{-\gamma}\,{\rm d}s\\
	&\le C t^{-\sigma(p)}\int_0^{\frac t2}(1+s)^{-\gamma}\,{\rm d}s+ C t^{-\gamma}\int_{\frac t2}^t (t-s)^{-\sigma(p)}\,{\rm d}s\\
	&\le  C t^{-\sigma(p)}\int_0^{\frac t2}(1+s)^{-\gamma}\,{\rm d}s+ Ct^{-\sigma(p)+1-\gamma},
	\end{aligned}
\end{equation}
which yields the desired bound in exterior regions, since
\begin{equation}\label{eq:int.gamma}
\int_0^{\frac t2}(1+s)^{-\gamma}\,{\rm d}s\asymp\begin{cases}
	 t^{1-\gamma},&\gamma<1,\\
	 \log t,&\gamma=1,\\
	 1,&\gamma>1,
	 	\end{cases}
\end{equation}

Let us consider then the critical and supercritical cases, which only occur if $N\ge 2\beta$.
We make the decomposition $|u|\le {\rm I}+{\rm II}$, where
\[
\begin{aligned}
{\rm I}(x,t)&=\int_0^t\int_{\{|y|<\frac{|x|}2\}}Y(x-y,t-s)|f(y,s)|\,{\rm d}y{\rm d}s,\\
{\rm II}(x,t)&=\int_0^t\int_{\{|y|>\frac{|x|}2\}}Y(x-y,t-s)|f(y,s)|\,{\rm d}y{\rm d}s.
\end{aligned}
\]
If $|y|<|x|/2$, then  $|x-y|>|x|/2$. Thus, if moreover $|x|\ge\nu t^\theta$, then $|x-y|(t-s)^{-\theta}>\nu /2$. Hence, the exterior  bound \eqref{eq:exterior.estimate.Y} together with the size hypothesis~\eqref{eq-hypothesis f} yield
$$
\begin{aligned}
{\rm I}(x,t)&\le  C \int_0^t (t-s)^{2\a-1}\int_{\{|y|<\frac{|x|}2\}}|x-y|^{-(N+2\beta)}|f(y,s)|\,{\rm d}y{\rm d}s\\
&\le  C |x|^{-(N+2\beta)}\int_0^t (t-s)^{2\a-1}\int_{\{|y|<\frac{|x|}2\}}|f(y,s)|\,{\rm d}y{\rm d}s\\
&\le C|x|^{-(N+2\beta)}\int_0^t (t-s)^{2\a-1}(1+s)^{-\gamma}\,{\rm d}s,
\end{aligned}
$$
and therefore
\[
\begin{aligned}
\|{\rm I}(\cdot,t)\|_{L^p(\{|x|\ge\nu t^\theta\})}
&\le C t^{-\sigma(p)}\int_0^{\frac t2}(1+s)^{-\gamma}\,{\rm d}s+C t^{-\sigma(p)+1-\gamma-2\a}\int_{\frac t2}^t(t-s)^{2\a-1}\,{\rm d}s\\
&\le C t^{-\sigma(p)}\int_0^{\frac t2}(1+s)^{-\gamma}\,{\rm d}s+C t^{-\sigma(p)+1-\gamma},
\end{aligned}	
\]
from where the desired bound for ${\rm I}$ follows using~\eqref{eq:int.gamma}.

In order to estimate ${\rm II}$ we use the decay condition~\eqref{eq:decay.condition} on $f$ and the integrability estimate~\eqref{eq:p.norm.Y} for the kernel with $p=1$. Since $\sigma(1)=1-\alpha$,
$$
\begin{aligned}
{\rm II}(x,t)&\le C\int_0^t(1+s)^{-\gamma}\int_{\{|y|>\frac{|x|}2\}}Y(x-y,t-s)|y|^{-N}\,{\rm d}y{\rm d}s\\
&\le C|x|^{-N}\int_0^t(1+s)^{-\gamma}\int_{\mathbb{R}^N}Y(x-y,t-s)\,{\rm d}y{\rm d}s= C|x|^{-N}\int_0^t(1+s)^{-\gamma}(t-s)^{\a-1}\,{\rm d}s,
\end{aligned}
$$
and therefore, remembering that $p>1$ in the critical and supercritical cases,
\[
\begin{aligned}
\|{\rm II}(\cdot,t)\|_{L^p(\{|x|\ge\nu t^\theta\})}
&\le C t^{-\sigma(p)}\int_0^{\frac t2}(1+s)^{-\gamma}\,{\rm d}s+C t^{-\sigma(p)+1-\gamma-\a}\int_{\frac t2}^t(t-s)^{\a-1}\,{\rm d}s\\
&\le C t^{-\sigma(p)}\int_0^{\frac t2}(1+s)^{-\gamma}\,{\rm d}s+C t^{-\sigma(p)+1-\gamma},
\end{aligned}	
\]
from where the desired bound for ${\rm II}$ follows using~\eqref{eq:int.gamma}.

\medskip

\noindent\textsc{Optimality. } We choose $f(x,t)=(1+t)^{-\gamma}\chi_{B_1}(x)$.

Let $t\ge1$, $s\in(0,t/2)$, $|x|<\mu t^\theta$, with $\mu>\nu$, and $|y|<1$. Then
\[
\frac{|x-y|}{(t-s)^\theta}\le\frac{1+\mu t^\theta}{(t/2)^\theta}
\le C,
\]
and hence $G\big((x-y)(t-s)^{-\theta}\big)\ge c>0$, since $G$ is positive. Therefore,
		\[
		u(x,t)\ge C\int_0^{\frac t2}(1+s)^{-\gamma}(t-s)^{-\sigma_*}\,{\rm d}s\ge Ct^{-\sigma_*}\int_0^{\frac t2}(1+s)^{-\gamma}\,{\rm d}s. \]
Thus,
\[
\begin{aligned}
		\|u(\cdot,t)\|_{L^p(\{|x|>\nu t^\theta\})}&\ge  \|u(\cdot,t)\|_{L^p(\{\mu >|x|/t^\theta>\nu \})}\ge C|\{\nu <|x|/t^\theta<\mu \}\big|^{1/p}t^{-\sigma_*}\int_0^{\frac t2}(1+s)^{-\gamma}\,{\rm d}s\\
		&=Ct^{-\sigma(p)}\int_0^{\frac t2}(1+s)^{-\gamma}\,{\rm d}s,
\end{aligned}
\]
which implies the desired lower bound.	
\end{proof}
\section{Compact regions}\label{sect-compacts}
 \setcounter{equation}{0}

In this section we study the decay/growth rate of the $L^p$ norm of the mild solution to \eqref{eq-hypothesis f} in compact sets.
\begin{proof}[Proof of Theorem \ref{teo-compactos}]
We have the bound $|u|\le {\rm I}+{\rm II}$, where
 	\[\begin{aligned}
 	{\rm I}(x,t)&=\int_0^{t-1}\int Y(x-y,t-s)|f(y,s)|\,{\rm d}y{\rm d}s,\\
 	{\rm II}(x,t)&=\int_{t-1}^t\int Y(x-y,t-s)|f(y,s)|\,{\rm d}y{\rm d}s.
\end{aligned}\]	

In order to estimate ${\rm I}$ we have to consider two cases, to take into account the different behaviors of the profile $G$ at the origin, depending on the dimension.	

If $1\le N<4\beta$, $G\in L^\infty(\R^N)$. Hence, using the decay assumption~\eqref{eq-hypothesis f}, for all $t\ge2$ we have
	 \begin{equation*}\label{eq-bound compact I1}\begin{aligned}
	 {\rm I}(x,t)&\le C\int_0^{t-1}(t-s)^{-\sigma_*}(1+s)^{-\gamma}\,{\rm d}s\le C t^{-\sigma_*}\int_0^{\frac t2}(1+s)^{-\gamma}\,{\rm d}s+ C t^{-\gamma}\int_{\frac t2}^{t-1}(t-s)^{-\sigma_*}\,{\rm d}s.
\end{aligned}
\end{equation*}
Using~\eqref{eq:int.gamma} and
\begin{equation}
\label{eq:int.dimension}
\int_{\frac t2}^{t-1}(t-s)^{-\sigma_*}\,{\rm d}s=\int_{1}^{\frac t2}\tau^{-\sigma_*}\,{\rm d}\tau\asymp\begin{cases}
	 t^{-\sigma_*+1},&1\le N<2\beta,\\
	 \log t,&N=2\beta,\\
	 1,&2\beta< N<4\beta,
	 	\end{cases}
\end{equation}
we get the desired bound for ${\rm I}$ for $p=\infty$, and therefore for all values of $p$, because we are considering compact sets.

If $N=4\beta$, we make the decomposition ${\rm I}= {\rm I}_1+{\rm I}_2+{\rm I}_3$, where
	 	\[\begin{aligned}
	 	{\rm I}_1(x,t)&= \int_0^{t-1}\int_{\{|x-y|<1\}} Y(x-y,t-s)|f(y,s)|\,{\rm d}y{\rm d}s,\\
	 	{\rm I}_2(x,t)&= \int_0^{t-1}\int_{\{1<|x-y|<(t-s)^\theta\}} Y(x-y,t-s)|f(y,s)|\,{\rm d}y{\rm d}s,
\\{\rm I}_3(x,t)&= \int_0^{t-1}\int_{\{|x-y|>(t-s)^\theta\}} Y(x-y,t-s)|f(y,s)|\,{\rm d}y{\rm d}s.
	 	\end{aligned}
	 	\]
Let $q=1$ if $p\in[1,p_{\rm c})$, $q>q_{\rm c}(p)$ as in~\eqref{eq:q.integrability.f} if $p\ge p_{\rm c}$.  Let $r$ satisfy $1+\frac1p=\frac1q+\frac1r$. Then $r\in [1,p_{\rm c})$ and,
	 	using the  behavior at the origin~\eqref{eq:estimates.G.origin} of the profile $G$, if~$s\le t-1$ we have
$$
\begin{aligned}
	\|Y(\cdot,t-s)\|_{L^r(B_1)}&\le C(t-s)^{-\sigma_*}\Big(\int_{B_1}\big(1-\log |x|+\log(t-s)\big)^r\,{\rm d}x\Big)^{1/r}\\
	&\le C(t-s)^{-\sigma_*}(1+\log(t-s)).
\end{aligned}
$$
Therefore, the integral condition~\eqref{eq:q.integrability.f} implies, for all $t\ge 2$,
$$
\begin{aligned}
\|{\rm I}_1(\cdot,t)\|_{L^p(\R^N)}&\le C\int_0^{t-1}\|Y(\cdot,t-s)\|_{L^r(B_1)}\|f(\cdot,s)\|_{L^q(\mathbb{R}^N)}\,{\rm d}s\\
&\le Ct^{-\sigma_*}\log t\int_0^{\frac t2}(1+s)^{-\gamma}\,{\rm d}s+C t^{-\gamma}\int_{\frac t2}^{t-1}(t-s)^{-\sigma_*}(1+\log(t-s))\,{\rm d}s.
\end{aligned}
$$
Since $\sigma_*=1+\alpha$ in this case, $\displaystyle\int_{\frac t2}^{t-1}(t-s)^{-\sigma_*}(1+\log(t-s))\,{\rm d}s\le C$ for all $t\ge 2$, and using also~\eqref{eq:int.gamma}, we finally get
$$
\|{\rm I}_1(\cdot,t)\|_{L^p(\R^N)}\le C\begin{cases}
	 	t^{-\gamma},&\gamma<\sigma_*,\\
	 	t^{-\sigma_*}\log t,&\gamma\ge\sigma_*.
	 	\end{cases}
$$

In order to bound ${\rm I}_2$, we observe that   $(t-s)^{-\theta}<|x-y|(t-s)^{-\theta}<1$ in the region under consideration. Hence, using the size condition~\eqref{eq-hypothesis f} and the behavior at the origin~\eqref{eq:estimates.G.origin} of the profile $G$,
	 	\[
	 	\begin{aligned}
	 	{\rm I}_2(x,t)&\le C\int_0^{t-1}(t-s)^{-\sigma_*}\int_{\{1<|x-y|<(t-s)^\theta\}}(1+\log(t-s))|f(y,s)|\,{\rm d}y{\rm d}s\\
	 	&\le Ct^{-\sigma_*}\log t\int_0^{\frac t2}(1+s)^{-\gamma}\,{\rm d}s+C t^{-\gamma}\int_{\frac t2}^{t-1}(t-s)^{-\sigma_*}(1+\log(t-s))\,{\rm d}s\\
&\le  C\begin{cases}
	 	t^{-\gamma},&\gamma<\sigma_*,\\
	 	t^{-\sigma_*}\log t,&\gamma\ge\sigma_*.
	 	\end{cases}
	 \end{aligned}
	 \]

	 	Finally, since $G$ is bounded outside the origin, using the size condition~\eqref{eq-hypothesis f},
	 	\[
	 	\begin{aligned}
	 	{\rm I}_3(x,t)
&\le C\int_0^{t-1}(t-s)^{-\sigma_*}(1+s)^{-\gamma}\,{\rm d}s\\
	 	&\le C t^{-\sigma_*}\int_0^{\frac t2}(1+s)^{-\gamma}\,{\rm d}s+Ct^{-\gamma}\int_{\frac t2}^{t-1}(t-s)^{-\sigma_*}\,{\rm d}s\le C 	 	t^{-\min\{\gamma,\sigma_*\}}.
\end{aligned}
\]

In order to bound ${\rm II}$, for all small dimensions, $1\le N\le 4\beta$, we take $r\in[1,p_{\rm c})$ as we did when we obtained the bound for ${\rm I}_1$. Notice that $\sigma(r)<1$. Hence, using~\eqref{eq:p.norm.Y} and~\eqref{eq:q.integrability.f} we get, for all $t\ge 2$,
\begin{equation*}
\label{eq-bound compactos II}
    \begin{aligned}
    \|{\rm II}(\cdot,t)\|_{L^p(\mathbb{R}^N)}
    &\le C\int_{t-1}^t\|Y(\cdot,t-s)\|_{L^r(\mathbb{R}^N)}\|f(\cdot,s)\|_{L^q(\R^N)}\,{\rm d}s\\
    &\le C\int_{t-1}^t(t-s)^{-\sigma(r)}(1+s)^{-\gamma}\,{\rm d}s\le C t^{-\gamma}\int_0^1 \tau^{-\sigma(r)}\,{\rm d}s=C t^{-\gamma},
    \end{aligned}	
\end{equation*}
which combined with the estimate for ${\rm I}$ yields the result.

\noindent\textsc{Optimality. }
We consider $f(x,t)=(1+t)^{-\gamma}\chi_{K+B_1}(x)$, where~$K$ is any compact set with measure different from 0.  If $x\in K$ and $|x-y|<1$, then $y\in K+B_1$. Therefore, for all $x\in K$,
\[
\begin{aligned}
u(x,t)&\ge \int_{0}^{t-1}\int_{\{|x-y|<1\}}Y(x-y,t-s) f(y,s)\,{\rm d}y{\rm d}s\\
&=\int_{0}^{t-1}(1+s)^{-\gamma}\int_{\{|x-y|<1\}}Y(x-y,t-s)\,{\rm d}y{\rm d}s.
\end{aligned}
\]
On the other hand, if $|x-y|<1$ and $s<t-1$, then
$|x-y|(t-s)^{-\theta}< 1$, and we have, using the behavior~\eqref{eq:estimates.G.origin} for $G$ close to the origin, that there exists a constant $c>0$ such that
$$
G\big((x-y)(t-s)^{-\theta}\big)\ge c
\begin{cases}
1,&1\le N<4\beta,\\
\log (t-s),& N=4\beta.
\end{cases}
$$
Hence, for all $x\in K$ and $t\ge2$,
\[
u(x,t)\ge c\begin{cases}\displaystyle
t^{-\sigma_*}\int_0^{\frac t2}(1+s)^{-\gamma}\,{\rm d}s+t^{-\gamma}\int_{\frac t2}^{t-1} (t-s)^{-\sigma_*}\,{\rm d}s,\quad&1\le N<4\beta,\\
\displaystyle t^{-\sigma_*}\log t\int_0^{\frac t2}(1+s)^{-\gamma}\,{\rm d}s+t^{-\gamma}\int_{\frac t2}^{t-1}(t-s)^{-\sigma_*}\log (t-s)\,{\rm d}s,&N=4\beta,
\end{cases}
\]
from where the desired lower bounds follow using~\eqref{eq:int.gamma} and~\eqref{eq:int.dimension},
since
$$
\int_{\frac t2}^{t-1}(t-s)^{-\sigma_*}\log (t-s)\,{\rm d}s\ge \int_1^2\tau^{-\sigma_*}\log \tau\,{\rm d}\tau=c>0
$$
for all $t\ge 4$ if $N=4\beta$.
\end{proof}
  \section{Intermediate scales}\label{sect-intermediate}
  \setcounter{equation}{0}
In this section we study the decay/growth rate of the $L^p$ norm of the mild solution to \eqref{eq-hypothesis f} in regions where $|x|\asymp g(t)$ with $g(t)\to\infty$ such that $g(t)=o(t^\theta)$.
\begin{proof}[Proof of Theorem \ref{teo-intermediate}]
We have $|u|\le {\rm I}+{\rm II}$, where
\begin{equation}
\label{eq-decomposition intermediate1}
\begin{aligned}
{\rm I}(x,t)&=\int_0^t\int_{\{|y|>\frac{|x|}2\}}|f(x-y,t-s)|Y(y,s)\,{\rm d}y{\rm d}s,\\
{\rm II}(x,t)&=\int_0^t\int_{\{|y|<\frac{|x|}2\}}|f(x-y,t-s)|Y(y,s)\,{\rm d}y{\rm d}s.
\end{aligned}
\end{equation}
We make the decomposition ${\rm I}={\rm I}_1+{\rm I}_2+{\rm I}_3$, where
\begin{equation}\label{eq-decomposition intermediate2}\begin{aligned}
{\rm I}_1(x,t)&= \int_0^{\big(\frac{|x|}2\big)^{1/\theta}}\int_{\{|y|>\frac{|x|}2\}}|f(x-y,t-s)|Y(y,s)\,{\rm d}y{\rm d}s,\\
{\rm I}_2(x,t)&= \int_{\big(\frac{|x|}2\big)^{1/\theta}}^{\frac t2}\int_{\{|y|>\frac{|x|}2\}}|f(x-y,t-s)|Y(y,s)\,{\rm d}y{\rm d}s,\\
{\rm I}_3(x,t)&= \int_{\frac t2}^t\int_{\{|y|>\frac{|x|}2\}}|f(x-y,t-s)|Y(y,s)\,{\rm d}y{\rm d}s.
\end{aligned}
\end{equation}

We start by estimating ${\rm I}_1$. In the region under consideration we have $|y|s^{-\theta}>1$. Hence we can use the estimate~\eqref{eq:exterior.estimate.Y} for $Y$ in exterior domains, which combined with the decay condition~\eqref{eq-hypothesis f} yields, for $|x|\le \delta t^\theta$ with $\delta<2$,
\begin{equation}\label{eq:estimate.I1.intermediate}
\begin{aligned}
{\rm I}_1(x,t)&\le C  |x|^{-(N+2\beta)}\int_0^{\big(\frac{|x|}{2}\big)^{1/\theta}} s^{2\a-1}(1+t-s)^{-\gamma}\,{\rm d}s\\
&\le C |x|^{-(N+2\beta)} t^{-\gamma}\int_0^{\big(\frac{|x|}{2}\big)^{1/\theta}} s^{2\a-1}\,{\rm d}s= C |x|^{(1-\sigma_*)/\theta} t^{-\gamma}.
\end{aligned}
\end{equation}
We conclude that $\|{\rm I}_1(\cdot,t)\|_{L^p(\{\nu<|x|/g(t)<\mu\})}\le  C t^{-\gamma}g(t)^{(1-\sigma(p))/\theta}$.

Let us now consider ${\rm I}_j$, $j=2,3$. If $N<4\beta$, the profile $G$ is bounded. Hence, using  the decay condition~\eqref{eq-hypothesis f}, if $|x|\le \delta t^\theta$ with $\delta<\min\{1,2^{1-\theta}\}$ we have
 \begin{equation}
 \label{eq:I2.I3.intermediate}
 \begin{aligned}
 {\rm I}_{2}(x,t)&
 \le Ct^{-\gamma}\int_{\big(\frac{|x|}2\big)^{1/\theta}}^{\frac t2}s^{-\sigma_*}\,{\rm d}s\le Ct^{-\gamma}
 \begin{cases}
 t^{-\sigma_*+1},&1\le N<2\beta,\\
\log(t^\theta/|x|), &N=2\beta,\\
 |x|^{(1-\sigma_*)/\theta},& 2\beta<N< 4\beta,
 \end{cases}
\\
  {\rm I}_{3}(x,t)&\le C t^{-\sigma_*}\int_{\frac t2}^t(1+t-s)^{-\gamma}\,{\rm d}s.
\end{aligned}
\end{equation}
Therefore, using also~\eqref{eq:int.dimension},
$$	
\begin{aligned}
  \|{\rm I}_{2}(\cdot,t)\|_{L^p(\{\nu<|x|/g(t)<\mu\}))}&\le   C g(t)^{N/p} \begin{cases}
  t^{-\sigma_*+1-\gamma},&1\le N<2\beta,\\
 \log(t^\theta/g(t))t^{-\gamma},&N=2\beta,\\
 g(t)^{(1-\sigma_*)/\theta}t^{-\gamma},& 2\beta<N< 4\beta,
 \end{cases}
 \\
  \|{\rm I}_{3}(\cdot,t)\|_{L^p(\{\nu<|x|/g(t)<\mu\}))}&
\le  C
   g(t)^{N/p}\begin{cases}
  t^{-\sigma_*+1-\gamma},&\gamma<1,\\
  t^{-\sigma_*}\log t,&\gamma=1,\\
  t^{-\sigma_*},&\gamma>1.
  \end{cases}
\end{aligned}
$$

When $N=4\beta$,  using the inner and outer behaviors~\eqref{eq:estimates.G.origin}--\eqref{eq:exterior.estimate.G} of the profile $G$ (the latter one implying boundedness), and the decay condition~\eqref{eq-hypothesis f}, and remembering that $\sigma_*=1+\alpha>1$ in this case, if $|x|\le \delta t^\theta$ with $\delta<\min\{1,2^{1-\theta}\}$,
\begin{equation}
\label{eq:I2.I3.4beta}
\begin{aligned}
 {\rm I}_{2}(x,t)\le&\, C \int_{\big(\frac{|x|}2\big)^{1/\theta}}^{\frac t2}\int_{\{\frac{|x|}2<|y|<s^\theta\}}|f(x-y,t-s)|s^{-\sigma_*}\big(1+\log(s^{\theta}/|y|)\big)\,{\rm d}y{\rm d}s\\
  &+C \int_{\big(\frac{|x|}2\big)^{1/\theta}}^{\frac t2}\int_{\{|y|>s^\theta\}}|f(x-y,t-s)|s^{-\sigma_*}\,{\rm d}y{\rm d}s\\
\le&\, C \int_{\big(\frac{|x|}2\big)^{1/\theta}}^{\frac t2}(1+t-s)^{-\gamma}s^{-\sigma_*}\big(1+\log(2s^{\theta}/|x|)\big)\,{\rm d}s
 \\
\le&\, C t^{-\gamma}|x|^{(1-\sigma_*)/\theta}\int_1^{\frac{2^{1-\theta}t^\theta}{|x|}}\tau^{-(1+\frac{\sigma_*-1}\theta)}(1+\log \tau)\,{\rm d}\tau
\le C t^{-\gamma}|x|^{(1-\sigma_*)/\theta},
\end{aligned}
\end{equation}
\begin{equation}
\label{eq:I2.I3.4beta.2}
\begin{aligned}
{\rm I}_{3}(x,t)\le&\, C \int_{\frac t2}^t\int_{\{\frac{|x|}2<|y|<s^\theta\}}|f(x-y,t-s)|s^{-\sigma_*}\big(1+\log(s^{\theta}/|y|)\big)\,{\rm d}y{\rm d}s\\
&+C \int_{\frac t2}^t\int_{\{|y|>s^\theta\}}|f(x-y,t-s)|s^{-\sigma_*}\,{\rm d}y{\rm d}s\\
 \le&\, C\big(1+\log(t^\theta/|x|)\big)t^{-\sigma_*}\int_{\frac t2}^t(1+t-s)^{-\gamma}\,{\rm d}s.
\end{aligned}
\end{equation}
Hence, using~\eqref{eq:int.dimension},
 $$
 \begin{aligned}
 \|{\rm I}_{2}(\cdot,t)\|_{L^p(\{\nu<|x|/g(t)<\mu \})}&\le   C t^{-\gamma}g(t)^{(1-\sigma(p))/\theta}\\
 \|{\rm I}_{3}(\cdot,t)\|_{L^p(\{\nu<|x|/g(t)<\mu \})}&\le   C g(t)^{N/p}\log(t^\theta/g(t))\begin{cases}
 t^{-\sigma_*+1-\gamma},&\gamma<1,\\
 t^{-\sigma_*}\log t,&\gamma=1,\\
 t^{-\sigma_*},&\gamma>1.
 \end{cases}
 \end{aligned}
$$
  	
Let us now turn to ${\rm II}$. We decompose it as ${\rm II}={\rm II}_1+{\rm II}_2$, where
$$
\begin{aligned}
{\rm II}_1(x,t)&=\int_0^{\frac t2}\int_{\{|y|<\frac{|x|}2\}}|f(x-y,t-s)|Y(y,s)\,{\rm d}y{\rm d}s,\\
{\rm II}_2(x,t)&=\int_{\frac t2}^t\int_{\{|y|<\frac{|x|}2\}}|f(x-y,t-s)|Y(y,s)\,{\rm d}y{\rm d}s.
\end{aligned}
$$

We start with ${\rm II}_1$. If $N<2\beta$, then  $G$ is bounded and $\sigma_*<1$. Hence, using the decay assumption~\eqref{eq-hypothesis f},
	\[\begin{aligned}
{\rm II}_{1}(x,t)&\le C	\int_0^{\frac t2}\int_{\{|y|<\frac{|x|}2\}}
	|f(x-y,t-s)|s^{-\sigma_*}\,{\rm d}y{\rm d}s\\
&\le C \int_0^{\frac t2}(1+t-s)^{-\gamma}s^{-\sigma_*}\,{\rm d}s\le C t^{-\gamma} \int_0^{\frac t2}s^{-\sigma_*}\,{\rm d}s= C t^{-\sigma_*+1-\gamma},
	\end{aligned}\]	
so that $\|{\rm II}_{1}(\cdot,t)\|_{L^p(\{\nu <|x|/g(t)<\mu\})}\le C g(t)^{N/p}t^{-\sigma_*+1-\gamma}$.
	
For  $2\beta\le N\le 4\beta$ we make the decomposition ${\rm II}_1={\rm II}_{11}+{\rm II}_{12}$, where
\[
\begin{aligned}
{\rm II}_{11}(x,t)&=\int_0^{\big(\frac{|x|}2\big)^{1/\theta}}\int_{\{|y|<\frac{|x|}{2}
\}}|f(x-y,t-s)|Y(y,s)\,{\rm d}y{\rm d}s,\\
{\rm II}_{12}(x,t)&=
\int_{\big(\frac{|x|}2\big)^{1/\theta}}^{\frac t2}\int_{\{|y|<\frac{|x|}{2}\}}|f(x-y,t-s)|Y(y,s)\,{\rm d}y{\rm d}s.
\end{aligned}
\]

If $p$ is subcritical, using the decay assumption~\eqref{eq-hypothesis f} we get (remember that $\sigma(p)<1$ in this case)
$$
\begin{aligned}
\|{\rm II}_{11}(\cdot,t)\|_{L^p(\{\nu<|x|/g(t)<\mu\})}&\le
C\int_0^{(\frac\mu2g(t))^{1/\theta}}\|f(\cdot,t-s)\|_{L^1(\R^N)}s^{-\sigma_*}\Big(\int_{\mathbb{R}^N}G^p(ys^{-\theta})\,{\rm d}y\Big)^{1/p}\,{\rm d}s\\
&\le C \int_0^{(\frac\mu2g(t))^{1/\theta}}(1+t-s)^{-\gamma}s^{-\sigma(p)}\,{\rm d}s\le C t^{-\gamma}g(t)^{(1-\sigma(p))/\theta},
\end{aligned}
$$
since $G\in L^p(\mathbb{R}^N)$ in this case.

If $p$ is not subcritical, since $|x-y|>|x|/2$ when $|y|<|x|/2$, using now the decay assumption~\eqref{eq:decay.condition}  we get,  for $|x|\le \delta t^\theta$ with $\delta<2^{1-\theta}$,
\begin{equation}
\label{eq:II11.intermediate}
\begin{aligned}
{\rm II}_{11}(x,t)&\le C|x|^{-N}\int_0^{(\frac{|x|}2)^{1/\theta}}\int_{\mathbb{R}^N}(1+t-s)^{-\gamma}s^{-\sigma_*}G(ys^{-\theta})\,{\rm d}y{\rm d}s
\\
&\le C|x|^{-N}t^{-\gamma}\int_0^{(\frac{|x|}2)^{1/\theta}}s^{-\sigma_*+N\theta}\,{\rm d}s\le C t^{-\gamma}|x|^{(1-\sigma_*)/\theta},
\end{aligned}
\end{equation}
since $G\in L^1(\mathbb{R}^N)$, and we get the same $L^p$ estimate as in the subcritical case.

As for ${\rm II}_{12}$, since $G$ is bounded if $N<4\beta$, using  the decay assumption~\eqref{eq-hypothesis f},
\begin{equation}
\label{eq:II12.N<4beta.intermediate}
{\rm II}_{12}(x,t)\le C\int_{(\frac{|x|}2)^{1/\theta}}^{\frac t2}(1+t-s)^{-\gamma}
	s^{-\sigma_*}\,{\rm d}s\le C\begin{cases}\displaystyle t^{-\gamma}\log(t^\theta/|x|),&N=2\beta,\\
\displaystyle t^{-\gamma}|x|^{(1-\sigma_*)/\theta},&2\beta< N<4\beta,
\end{cases}
\end{equation}
if $|x|\le \delta t^\theta$ with $\delta<2^{1-\theta}$, so that
\begin{equation*}
	\|{\rm II}_{12}(\cdot,t)\|_{L^p(\{\nu<|x|/g(t)<\mu\})}\le C\begin{cases}
t^{-\gamma} g(t)^{N/p}
\log(t^\theta/g(t)),&
N=2\beta,\\
t^{-\gamma} g(t)^{(1-\sigma(p))/\theta},&2\beta< N<4\beta.
\end{cases}
\end{equation*}

On the other hand, if $N=4\beta$,  since  $|y|\le s^\theta$ in the region under consideration,  using the behavior~\eqref{eq:estimates.G.origin} of $G$ close to the origin,
\begin{equation}
\label{eq:II12.intermediate}
{\rm II}_{12}(\cdot,t)\le
C\int_{(\frac\nu2g(t))^{1/\theta}}^{\frac t2}|f(x-y,t-s)|s^{-\sigma_*}\int_{\{|y|<\frac{|x|}2\}}\big(1+\big|\log(|y|/s^\theta)\big|\big)\,{\rm d}y{\rm d}s.
\end{equation}
It is easy to check that for each $p\in [1,\infty)$ there exists a constant $C$ such that
\begin{equation}\label{eq:estimate.integral.log}
\int_{\{|z|<\rho\}}\big(1+\big|\log|z|\big|\big)^p)\,{\rm d}y\le C\rho^N\big(1+|\log\rho|\big)^p\quad \text{for every } \rho\ge0.
\end{equation}
Therefore, using the decay assumption~\eqref{eq-hypothesis f}, for all
$p\in[1,p_{\rm c})$ we have
$$
\begin{aligned}
	\|{\rm II}_{12}(\cdot&,t)\|_{L^p(\{\nu<|x|/g(t)<\mu\})}\\
&\le
C\int_{(\frac\nu2g(t))^{1/\theta}}^{\frac t2}\|f(\cdot,t-s)\|_{L^1(\R^N)}s^{-\sigma_*}\Big(\int_{\{|y|<\frac\mu2g(t)\}}\big(1+\big|\log(|y|/s^\theta)\big|\big)^p\,{\rm d}y\Big)^{1/p}\,{\rm d}s\\
&\le Ct^{-\gamma}\int_{(\frac\nu2g(t))^{1/\theta}}^{\frac t2}s^{-\sigma(p)}\Big(\int_{\{|z|<\frac{\mu g(t)}{2s^\theta}\}}\big(1+\big|\log|z|\big|\big)^p)\,{\rm d}y\Big)^{1/p}\,{\rm d}s\\
	&\le C t^{-\gamma}g(t)^{N/p}\int_{(\frac\nu2g(t))^{1/\theta}}^{\frac t2}s^{-\sigma_*}\big(1+\big|\log(2s^\theta/(\mu g(t)))\big|\big)\,{\rm d}s\\
&= C t^{-\gamma}g(t)^{(1-\sigma(p))/\theta}\int_{\nu/\mu}^{\frac{2^{1-\theta}t^\theta}{\mu g(t)}}\tau^{-\big(1+(\sigma_*-1)/\theta\big)}|\log \tau|\,{\rm d}\tau\le C t^{-\gamma}g(t)^{(1-\sigma(p))/\theta}
\end{aligned}
$$
since $\sigma_*=1+\alpha>1$ in this case.

If $N=4\beta$ and $p\ge p_{\rm c}$, we plug the decay assumption~\eqref{eq:decay.condition} in~\eqref{eq:II12.intermediate}. Therefore, since $|x-y|>|x|/2$ if $|y|<|x|/2$, using~\eqref{eq:estimate.integral.log} with $p=1$,
\begin{equation}
\label{eq:II12.4=Nbeta.intermediate}
\begin{aligned}
{\rm II}_{12}(x,t)&\le
		C|x|^{-N}\int_{(\frac{|x|}2)^{1/\theta}}^{\frac t2}(1+t-s)^{-\gamma}s^{-\sigma_*}
		\int_{\{|y|<\frac{|x|}2\}} \big(1+\log(s^\theta/|y|)\big)\,{\rm d}y{\rm d}s\\
&\le 		Ct^{-\gamma}|x|^{-N}\int_{(\frac{|x|}2)^{1/\theta}}^{\frac t2}s^{-\sigma_*+N\theta}
		\int_{\{|z|<\frac{|x|}{2s^\theta}\}} \big(1+\big|\log|z|\big|\big)\,{\rm d}z{\rm d}s\\
	&\le Ct^{-\gamma}\int_{(\frac{|x|}2)^{1/\theta}}^{\frac t2}s^{-\sigma_*}\big(1+\log(2s^\theta/|x|)\big)\,{\rm d}s\\
		&	\le C t^{-\gamma} |x|^{(1-\sigma_*)/\theta}\int_1^{\frac{2^{1-\theta}t^\theta}{|x|}}\tau^{-\big(1+(\sigma_*-1)/\theta\big)}\log \tau\,{\rm d}\tau\le C t^{-\gamma} |x|^{(1-\sigma_*)/\theta},
	\end{aligned}
\end{equation}
and we get the same $L^p$ estimate as in the subcritical case.

 We finally consider ${\rm II}_2$. If $N<4\beta$, using that $G$ is bounded and the decay assumption~\eqref{eq-hypothesis f},
\begin{equation}
\label{eq:II2.N<4beta.intermediate}
\begin{aligned}
  {\rm II}_{2}(x,t)&\le C\int_{\frac t2}^t\int_{\{|y|<\frac{|x|}2\}}|f(x-y,t-s)|s^{-\sigma_*}\,{\rm d}y{\rm d}s\le C t^{-\sigma_*}\int_{\frac t2}^t(1+t-s)^{-\gamma}\,{\rm d}s,
  	\end{aligned}
\end{equation}
so that, using also~\eqref{eq:int.gamma}
\[	\|{\rm II}_{2}(\cdot,t)\|_{L^p(\{\nu<|x|/g(t)<\mu\})}\le
C g(t)^{N/p}\begin{cases}
t^{-\sigma_*+1-\gamma},&\gamma<1,\\
t^{-\sigma_*}\log t,&\gamma=1,\\
t^{-\sigma_*},&\gamma>1.
\end{cases}
\]
  	
Let now $N=4\beta$. Since $|y|<s^\theta$ in the region under consideration, using the behavior~\eqref{eq:estimates.G.origin} of $G$ close to the origin,
\begin{equation}
\label{eq:II2.N=4beta.intermediate}
  {\rm II}_{2}(x,t)\le C\int_{\frac t2}^t\int_{\{|y|<\frac{|x|}2\}}|f(x-y,t-s)|s^{-\sigma_*}\big(1+\big|\log(|y|/s^\theta)\big|\big)\,{\rm d}y{\rm d}s.
\end{equation}
Thus, using the decay assumption~\eqref{eq-hypothesis f} and the estimates~\eqref{eq:int.gamma} and~\eqref{eq:estimate.integral.log}, for $p\in [1,\infty)$ we have
 \[	\begin{aligned}
 \|{\rm II}_{2}(\cdot,t)\|_{L^p(\{\nu<|x|/g(t)<\mu\})}&\le
 C\int_{\frac t2}^t\|f(\cdot,t-s)\|_{L^1(\R^N)}s^{-\sigma_*}\Big(\int_{\{|y|<\frac\mu2g(t)\}}\big(1+\big|\log(|y|/s^\theta)\big|\big)^p\,{\rm d}y\Big)^{1/p}\,{\rm d}s   \\
 &\le C g(t)^{N/p}\log\big(t^\theta/g(t)\big) t^{-\sigma_*}\int_{\frac t2}^t(1+t-s)^{-\gamma}\,{\rm d}s\\
 &\le C g(t)^{N/p}\log\big(t^\theta/g(t)\big)\begin{cases}
 t^{-\sigma_*+1-\gamma},&\gamma<1,\\
 t^{-\sigma_*}\log t,&\gamma=1,\\
 t^{-\sigma_*},&\gamma>1.
 \end{cases}
 \end{aligned}\]

 If $N=4\beta$ and $p=\infty$, we use the decay assumption~\eqref{eq:decay.condition} instead, so that, as $|x-y|>|x|/2$ in the region under consideration, if $|x|\le t^\theta$ we have
\begin{equation}
\label{eq:II2.N=4beta.pinfty}
\begin{aligned}
 {\rm II}_{2}(x,t)&\le C|x|^{-N}\int_{\frac t2}^t\int_{\{|y|<\frac{|x|}2\}}(1+t-s)^{-\gamma} s^{-\sigma_*}\big(1+\big|\log (|y|/s^\theta)\big|\big)\,{\rm d}y{\rm d}s\\
 &\le Ct^{-\sigma_*+N\theta}|x|^{-N}\int_{\frac t2}^t\int_{\{|z|<\frac{|x|}{2 s^\theta}\}}(1+t-s)^{-\gamma} \big(1+\big|\log |z|\big|\big)\,{\rm d}z{\rm d}s\\
 &\le Ct^{-\sigma_*} \log(t^\theta/|x|)\int_{\frac t2}^t(1+t-s)^{-\gamma}\,{\rm d}s.
 \end{aligned}
 \end{equation}	
Therefore, if $|x|\asymp g(t)$,
${\rm II}_2(x,t)
 \le C\log\big(t^\theta/g(t)\big)\begin{cases}
 t^{-\sigma_*+1-\gamma},&\gamma<1,\\
 t^{-\sigma_*}\log t,&\gamma=1,\\
 t^{-\sigma_*},&\gamma>1.
 \end{cases}
$

\medskip

\noindent\textsc{Optimality. }  To see that the estimates are sharp, we consider $f(x,t)=(1+t)^{-\gamma}\chi_{B_1}(x)$, so that
\[
u(x,t)= \int_0^t\int_{\{|x-y|<1\}}(1+t-s)^{-\gamma}Y(y,s)\,{\rm d}y{\rm d}s.
\]

If $|x-y|<1$, then $|y|\le |x|+1\le 2|x|$ if $t$ is large enough, since $|x|\ge \nu g(t)$ and $g(t)\to\infty$. In particular, if $s\ge t/2$, then $|y|/s^\theta\le2^{1+\theta}\mu g(t)/t^\theta\le 1$ for $t$ large, since $g(t)=o(t^\theta)$. Hence, using the behavior~\eqref{eq:estimates.G.origin} of $G$ close to the origin,
\[
\begin{aligned}
u(x,t)
&\ge c \begin{cases}
\displaystyle
\int_{\frac t2}^t (1+t-s)^{-\gamma}s^{-\sigma_*}\,{\rm d}s,& N<4\beta,\\[14pt]
\displaystyle\int_{\frac t2}^t\int_{\{|x-y|<1\}}(1+t-s)^{-\gamma}s^{-\sigma_*}\big(1+\log(s^\theta/|y|)\big)\,{\rm d}y{\rm d}s,&N=4\beta,
\end{cases}
\\
&\ge c t^{-\sigma_*}\int_{\frac t2}^t (1+t-s)^{-\gamma}\,{\rm d}s
\begin{cases}
1,& N<4\beta,\\
\log(t^\theta/|x|),&N=4\beta,
\end{cases}
\end{aligned}\]
so that
\[
\|u(\cdot,t)\|_{L^p(\{\nu<|x|/g(t)<\mu\})}\ge c g(t)^{N/p}
\begin{cases}
t^{-\sigma_*+1-\gamma},&\gamma<1,\\
t^{-\sigma_*}\log t,&\gamma=1,\hskip2cm N<4\beta,\\
t^{-\sigma_*},&\gamma>1,\\[0.3cm]
\log(t^\theta/g(t))t^{-\sigma_*+1-\gamma},&\gamma<1,\\
\log(t^\theta/g(t))t^{-\sigma_*}\log t,&\gamma=1,\hskip2cm N=4\beta,\\
\log(t^\theta/g(t))t^{-\sigma_*},&\gamma>1.
\end{cases}\]

On the other hand, since $G$ is positive, for $t$ large so that $|x|>1$,
\begin{equation}
\label{eq:estimate.below.intermediate}
\begin{aligned}
u(x,t)&\ge c\int_{|x|^{1/\theta}}^{t}(1+t-s)^{-\gamma}s^{-\sigma_*}\,{\rm d}s\ge ct^{-\gamma}\int_{|x|^{1/\theta}}^ts^{-\sigma_*}\,{\rm d}s\\
&=c t^{-\gamma}\begin{cases}
\log\big(t^\theta/|x|\big),& N=2\beta,\\
|x|^{(1-\sigma_*)/\theta},&2\beta<N\le 4\beta.
\end{cases}
\end{aligned}
\end{equation}
Therefore, as $p_{\rm c}=\infty$ when $N=2\beta$,
\[
\|u(\cdot,t)\|_{L^p(\{\nu<|x|/g(t)<\mu\})}\ge c g(t)^{N/p}t^{-\gamma}
\begin{cases}
\log(t^\theta/g(t)),&N=2\beta,\\
g(t)^{(1-\sigma_*)/\theta},&2\beta<N\le 4\beta.
\end{cases}\]
Combining both estimates, we get the desired lower bounds.
\end{proof}

 \section{Global estimates}\label{sect-conexion}
  \setcounter{equation}{0}
In this section we establish the  behavior of the \emph{global} $L^p(\R^N)$ norms of the mild solution to~\eqref{eq-f}.
\begin{proof}[Proof of Theorem~\ref{teo-conexion}]    	
The behavior in the subcritical cases, with $p$ satisfying~\eqref{eq:subcritical.range}, was already established in the proof of Theorem~\ref{teo-exterior}; see estimate~\eqref{eq:global.subcritical}. Hence we turn our attention to the critical and supercritical cases.
Due to  the results of theorems \ref{teo-exterior} and \ref{teo-compactos}, it is enough to show that the estimates are true in some region of the form $\{R\le |x|\le \delta t^\theta\}$ with  $R,\delta>0$.

We will follow the computations in the proof of Theorem \ref {teo-intermediate}.

We have $|u|\le{\rm I}+{\rm II}$, with ${\rm I}$ and ${\rm II}$ as in~\eqref{eq-decomposition intermediate1}. The term ${\rm I}$ is further decomposed as ${\rm I}={\rm I}_1+{\rm I}_2+{\rm I}_3$, with ${\rm I}_j$, $j\in\{1,2,3\}$ as in~\eqref{eq-decomposition intermediate2}.
Since $|x|<\delta t^\theta$ in the region we are interested in, taking  $\delta<\min\{1,2^{1-\theta}\}$ we have~\eqref{eq:estimate.I1.intermediate}, and hence, since $(1-\sigma_*)/\theta=2\beta-N=-N/p_{\rm c}$,
\[	
\|{\rm I}_{1}(\cdot,t)\|_{L^p(\{R<|x|<\delta t^\theta\})}\le     C
\begin{cases}
t^{-\gamma}\log t,&p=p_{\rm c},\\
t^{-\gamma},&p>p_{\rm c}.
\end{cases}
\]

We consider now ${\rm I}_{2}$ and ${\rm I}_3$.
If $2\beta\le N<4\beta$ we have~\eqref{eq:I2.I3.intermediate}, from where we immediately obtain
\begin{equation}
\label{eq:estimate.I2.global}
\begin{aligned}
\|{\rm I}_{2}(\cdot,t)\|_{L^p(\{R<|x|<\delta t^\theta\})}&\le C
\begin{cases}
t^{-\gamma}\log t,&p=p_{\rm c},\\
t^{-\gamma},&p>p_{\rm c},
\end{cases}\\
\|{\rm I}_{3}(\cdot,t)\|_{L^p(\{R<|x|<\delta t^\theta\})}&\le   C
\begin{cases}
t^{-\sigma(p)+1-\gamma},&\gamma<1,\\
t^{-\sigma(p)}\log t,&\gamma=1,\\
t^{-\sigma(p)},&\gamma>1.
\end{cases}
\end{aligned}
\end{equation}
If
$N=4\beta$, we have~\eqref{eq:I2.I3.4beta}--\eqref{eq:I2.I3.4beta.2} instead, and we get  again~\eqref{eq:estimate.I2.global}.

Let us now turn to ${\rm II}$, always for values of $p$ which are not subcritical.  Using~\eqref{eq:II11.intermediate},~\eqref{eq:II12.N<4beta.intermediate}, and~\eqref{eq:II12.4=Nbeta.intermediate}, we obtain, for $j=1,2$,
$$
\|{\rm II}_{1j}(\cdot,t)\|_{L^p(\{R<|x|<\delta t^\theta\})}\le C
\begin{cases}
t^{-\gamma}\log t,&p=p_{\rm c},\\
t^{-\gamma},&p>p_{\rm c}.
\end{cases}
$$
On the other hand, if $N<4\beta$, using~\eqref{eq:II2.N<4beta.intermediate} we obtain
\[	
\|{\rm II}_{2}(\cdot,t)\|_{L^p(\{R<|x|<\delta t^\theta\})}\le
C \begin{cases}
t^{-\sigma(p)+1-\gamma},&\gamma<1,\\
t^{-\sigma(p)}\log t,&\gamma=1,\\
t^{-\sigma(p)},&\gamma=1.
\end{cases}
\]
If $N=4\beta$, starting from~\eqref{eq:II2.N=4beta.intermediate}, using the decay assumption~\eqref{eq-hypothesis f} and the estimate~\eqref{eq:int.gamma}, for $p\in [1,\infty)$ we have
 \[	\begin{aligned}
 \|{\rm II}_{2}(\cdot,t)\|_{L^p(\{R<|x|<\delta t^\theta\})}&\le
 C\int_{\frac t2}^t\|f(\cdot,t-s)\|_{L^1(\R^N)}s^{-\sigma_*}\Big(\int_{\{|y|<\delta t^\theta\}}\big(1+\big|\log(|y|/s^\theta)\big|\big)^p\,{\rm d}y\Big)^{1/p}\,{\rm d}s   \\
 &\le C  t^{-\sigma(p)}\int_{\frac t2}^t(1+t-s)^{-\gamma}\Big(\int_{\{|z|<\delta\}}\big(1+\big|\log|z|\big|\big)^p\,{\rm d}y\Big)^{1/p}\,{\rm d}s\\
 &\le C \begin{cases}
 t^{-\sigma(p)+1-\gamma},&\gamma<1,\\
 t^{-\sigma(p)}\log t,&\gamma=1,\\
 t^{-\sigma(p)},&\gamma>1.
 \end{cases}
 \end{aligned}\]
On the other hand, if $p=\infty$, for $R<|x|<\delta t^\theta$ estimate~\eqref{eq:II2.N=4beta.pinfty} yields
\[
{\rm II}_{2}(x,t)\le
 C \begin{cases}
   		t^{-\sigma_*+1-\gamma}\log t,&\gamma<1,\\
   		t^{-\sigma_*}(\log t)^2,&\gamma=1,\\
   		t^{-\sigma_*}\log t,&\gamma>1.
   		\end{cases}
\]

\medskip

\noindent\textsc{Optimality. } Except when we have simultaneously $p=p_{\rm c}$, $\gamma<1$, and  $N>2\beta$, the global rate coincides either with the one in exterior regions or the one in compact sets, which we have already shown to be optimal. To show that our estimates are sharp also in this exceptional case, we consider, once more, $f(x,t)=(1+t)^{-\gamma}\chi_{B_1}(x)$. The desired lower bound follows then immediately from~\eqref{eq:estimate.below.intermediate}.
\end{proof}

\medskip

\noindent\textbf{Acknowledgments. } This project has received funding from the European Union's Horizon 2020 research and innovation programme under the Marie Sklodowska-Curie grant agreement No.\,777822. Carmen Cort\'azar supported by  FONDECYT grant 1190102 (Chile). Fernando Quir\'os supported by Ministerio de Econom\'{\i}a y Competitividad (Spain), through project MTM2017-87596-P, and by Ministerio de Ciencia e Innovaci\'on (Spain), through project ICMAT-Centro de excelencia “Severo Ochoa”, CEX2019-000904-S. Noem\'{\i} Wolanski supported by CONICET PIP 11220150100032CO 2016-2019, UBACYT 20020150100154BA, ANPCyT
PICT2016-1022 and MathAmSud 13MATH03 (Argentina).


\begin{thebibliography}{99}
	

	
\bibitem{Caputo-1967} Caputo, M. \emph{Linear models of dissipation whose $Q$ is almost frequency independent--II.} Geophys.
	J. R. Astr. Soc. 13 (1967), 529--539.
	

\bibitem{Cartea-delCastilloNegrete-2007} Cartea, \'A.; del Castillo-Negrete, D. \emph{Fluid limit of the continuous-time random walk with general Lévy jump distribution functions.} Phys. Rev. E 76 (2007) 041105.

\bibitem{Compte-Caceres-1998} Compte, A.; C\'aceres, M.\,O. \emph{Fractional dynamics in random velocity fields.} Phys. Rev. Lett. 81 (1998) 3140--3143.

\bibitem{Cortazar-Quiros-Wolanski-2021a} Cortazar, C.; Quir\'os, F.; Wolanski, N. \emph{A heat equation with memory: large-time behavior.} J. Funct. Anal. 281 (2021), no.\,9, 109174.

\bibitem{Cortazar-Quiros-Wolanski-2021b} Cort\'azar, C.; Quir\'os, F.; Wolanski, N. \emph{Large-time behavior for a fully nonlocal heat equation.} Vietnam J. Math. 49 (2021), no.\,3, 831--844.

\bibitem{Cortazar-Quiros-Wolanski-2022} Cortazar, C.; Quir\'os, F.; Wolanski, N.     \emph{Decay/growth rates for inhomogeneous heat equations with memory. The case of large dimensions.}   Math. Eng. 4 (2022), no.\,3, 1--17.


\bibitem{Cortazar-Quiros-Wolanski-2021-Preprint} Cort\'azar, C.; Quir\'os, F.; Wolanski, N. \emph{Asymptotic profiles for inhomogeneous heat equations with memory}. Preprint.

\bibitem{delCastilloNegrete-Carreras-Lynch-2004} del Castillo-Negrete, D.; Carreras, B.\,A.; Lynch, V.\,E. \emph{Fractional diffusion in plasma turbulence.}
Physics of Plasmas 11 (2004), no.\,8,  3854--3864.

\bibitem{delCastilloNegrete-Carreras-Lynch-2005} del Castillo-Negrete, D.; Carreras, B.\,A.; Lynch, V.\,E. \emph{Nondiffusive transport in plasma turbulence:
A fractional diffusion approach.} Physical Review Letters 94 (2005), no.\,6,  065003.


\bibitem{Dzherbashyan-Nersesian-1968} Dzherbashyan, M.\,M.; Nersesian, A.\,B. \emph{Fractional derivatives and the Cauchy problem for
differential equations of fractional order}.  (Russian) Izv. Akad. Nauk Arm. SSR, Mat. 3 (1968), 3–-29.

	
	\bibitem{Eidelman-Kochubei-2004}
	Eidelman, S.\,D.; Kochubei, A.\,N. \emph{Cauchy problem for fractional diffusion equations.} J. Differential Equations 199 (2004), no.\,2, 211--255.
	


\bibitem{Gerasimov-1948}
Gerasimov, A.\,N. \emph{A generalization of linear laws of deformation and its application to problems of internal friction}. (Russian) Akad. Nauk SSSR. Prikl. Mat. Meh. 12 (1948), 251--260.





\bibitem{Gross-1947} Gross, B. \emph{On creep and relaxation}. J. Appl. Phys. 18 (1947), 212--221.


%
	
	
	
	\bibitem{Kemppainen-Siljander-Zacher-2017} Kemppainen, J.; Siljander, J.; Zacher, R. \emph{Representation of solutions and large-time behavior for fully nonlocal diffusion equations.} J. Differential Equations 263 (2017), no.\,1, 149--201.
	
%



\bibitem{Liouville-1832} Liouville, J. \emph{Memoire sur quelques questions de g\'eometrie et de m\'eecanique, et sur un nouveau
gentre pour resoudre ces questions}. J. Ecole Polytech. 13 (1832), 1--69.


	
	
	\bibitem{Metzler-Klafter-2000}	Metzler, R.; Klafter, J. \emph{The random walk's guide to anomalous diffusion: a fractional dynamics approach.} Phys. Rep. 339 (2000), no.\,1, 77 pp.

\bibitem{Pruss-book} Pr\"uss, J. \lq\lq Evolutionary integral equations and applications''. [2012] reprint of the 1993 edition. Modern Birkh\"auser Classics. Birkh\"auser/Springer Basel AG, Basel, 1993.  ISBN: 978-3-0348-0498-1.

\bibitem{Rabotnov-1966}    Rabotnov, Yu.\,N. \lq\lq Polzuchest Elementov Konstruktsii''. (Russian) Nauka, Moscow (1966); English translation:
\lq\lq Creep Problems in Structural Members''. North-Holland, Amsterdam (1969).


	

\bibitem{Stein-book-1970} Stein, E. M. \lq\lq Singular integrals and differentiability properties of functions''. Princeton Mathematical Series, No. 30 Princeton University Press, Princeton, N.J. 1970.
	

\bibitem{Zaslavsky-2002} Zaslavsky, G.\,M. \emph{Chaos, fractional kinetics, and anomalous transport.} Phys. Rep. 371 (2002), no.\,6, 461--580.
	
\end{thebibliography}
\end{document}